\documentclass[oneside,11pt]{article}
\usepackage{graphicx, subfigure}
\usepackage[top=1in,bottom=1in,left=1in,right=1in]{geometry}
\usepackage{amsfonts, amsmath, amssymb, amsthm, constants, bbm}
\usepackage{enumitem}

\usepackage[mathscr]{euscript}
\usepackage{pifont}
\usepackage{mathtools}
\mathtoolsset{showonlyrefs}
\usepackage{natbib}
\usepackage{mathrsfs}

\usepackage{multirow}
\usepackage{float}
\usepackage{color}
\definecolor{darkred}{RGB}{100,0,0}
\definecolor{darkgreen}{RGB}{0,100,0}
\definecolor{darkblue}{RGB}{0,0,150}
\usepackage{hyperref}
\hypersetup{colorlinks=true, linkcolor=darkred, citecolor=darkgreen, urlcolor=darkblue}

\def\beq{\begin{equation}} 
\def\eeq{\end{equation}}
\def\beqn{\begin{eqnarray*}}
\def\eeqn{\end{eqnarray*}}
\def\Bitem{\begin{itemize}\setlength{\itemsep}{.2in}}
\def\bitem{\begin{itemize}\setlength{\itemsep}{.05in}}
\def\eitem{\end{itemize}}
\def\Benum{\begin{enumerate}\setlength{\itemsep}{.2in}}
\def\benum{\begin{enumerate}\setlength{\itemsep}{.05in}}
\def\eenum{\end{enumerate}}
\def\bmult{\begin{multline*}}
\def\emult{\end{multline*}}
\def\bcenter{\begin{center}}
\def\ecenter{\end{center}}
\def\bframe{\begin{frame}}
\def\eframe{\end{frame}}

\def\cA{\mathcal{A}}

\def\cM{\mathcal{M}}

\def\bbG{\mathbb{G}}

\newcommand{\E}{\operatorname{\mathbb{E}}}
\renewcommand{\P}{\operatorname{\mathbb{P}}}

\def\eps{\varepsilon}

\def\1{\mathbbm{1}}

\def\sF{\mathscr{F}}
\def\sG{\mathscr{G}}

\newtheorem{thm}{Theorem}
\newtheorem{prp}{Proposition}

\newtheorem{assump}{Assumption}

\theoremstyle{remark}
\newtheorem{Def}{Definition}

\title{\textbf{Dynamic Treatment Effects under Functional Longitudinal Studies}\footnote{This preprint is a draft of my original ideas, with so many errors, and is being shared for the sole purpose
of allowing me to cite it as a source of my own work. I do not intend to publish this preprint in any formal
manner at this stage but retain all rights to my ideas. The preprint may be cited as a reference, but should
not be reproduced or distributed without my express permission. I take no responsibility for the accuracy
or reliability of this preprint and it should not be considered a final or official version of my work.}}

\author{Andrew Ying\\
Irvine, CA 92606, USA \\
\texttt{aying9339@gmail.com}
}
\date{}
\begin{document}

\maketitle
\begin{abstract}
Establishing causality is a fundamental goal in fields like medicine and social sciences. While randomized controlled trials are the gold standard for causal inference, they are not always feasible or ethical. Observational studies can serve as alternatives but introduce confounding biases, particularly in complex longitudinal data, where treatment-confounder feedback complicates analysis. The challenge increases with Dynamic Treatment Regimes (DTRs), where treatment allocation depends on rich historical patient data. The advent of real-time healthcare monitoring technologies, such as MIMIC-IV and Continuous Glucose Monitoring (CGM), has popularized Functional Longitudinal Data (FLD). However, there is yet no investigate of causal inference for FLD with DTRs. In this paper, we address it by developing a population-level framework for functional longitudinal data, accommodating DTRs. To that end, we define the potential outcomes and causal effects of interest. We then develop identification assumptions, and derive g-computation, inverse probability weighting, and doubly robust formulas through novel applications of stochastic process and measure theory. We further show that our framework is nonparametric and compute the efficient influence curve using semiparametric theory. Last, we illustrate our framework's potential through Monte Carlo simulations.
\end{abstract}



\section{Introduction}
Understanding treatment effects and establishing causality are fundamental goals in fields such as medicine and social sciences. While double-blinded randomized controlled trials are considered the gold standard for causal inference, they may not always be feasible or ethical. Observational studies are sometimes feasible alternative for comparative studies yet introduced additional confounding bias. Complex observational longitudinal data, pose additional challenges for causal inference due to the interplay among time-varying confounders, treatment, and outcome processes, also known as the treatment-confounder feedbacks \citep{robins1997marginal, robins1999association, hernan2020causal}. Another complication within longitudinal studies, due to the many possible treatment combinations across time and rich historical information of patient characteristics, treatment process, and prognoses, is the need to investigate causal effects under dynamic treatment regimes, that is, treatment allocation might depend on all historical information.

The advent of real-time monitoring technologies in healthcare has led to the continuous-time measurement of patient data. For example, the Medical Information Mart for Intensive Care IV (MIMIC-IV) \citep{johnson2023mimic} is a freely accessible electronic health record (EHR) database that records ICU care data, including physiological measurements, laboratory values, medication administration, and clinical events. Another example is Continuous Glucose Monitoring (CGM) \citep{rodbard2016continuous, klonoff2017continuous}, an increasingly adopted technology for insulin-requiring patients that provides insights into glycemic fluctuations. CGM offers a real-time, high-resolution stream of data, capturing the intricate fluctuations in interstitial fluid glucose levels every few minutes. These examples illustrate the recent prevalence of functional longitudinal data (FLD).

However, there is no investigation of causal inference for dynamic treatment regimes of longitudinal studies with functional data. This is possibly because of the fact that FLD, characterized by continuous-time processes and high-dimensional measurements, present additional complexities for causal inference due to the presence of an infinite number of potential treatment-confounder feedback loops and the lack of a joint density for these processes. Existing methods developed for classical longitudinal data can not be directly applicable to functional longitudinal data. For less distraction, we dedicate a thorough literature review later in Section \ref{sec:related}. Therefore, there is a pressing need for more sophisticated and nuanced methods of causal inference that can effectively handle the intricacies of functional longitudinal data.


We aim to fill this gap partially by developing a comprehensive population-level framework for causal inference in functional longitudinal data under dynamic treatment regimes.

To that end, we proceed with the following steps routinely adopted by classical causal inference frameworks for non-functional longitudinal data:
\begin{enumerate}
    \item We generalize the potential outcome framework \citep{neyman1923applications, rubin1974estimating, holland1986statistics} and define the causal parameter of interest as a functional of measures on unobserved potential outcomes, for FLD under general DTRs. Here by ``general'' we mean the DTRs can be either deterministic or stochastic, and be prespecified or allowed to depend on the actual observed TRs;
    \item We propose a set of identification assumptions tailored to the FLD setting, to identify the causal parameter of interest defined earlier. 
    Identification methods include a g-computation formula \citep{greenland1986identifiability, robins2000robust, bang2005doubly}, an inverse probability weighting (IPW) formula \citet{rosenbaum1983central, hernan2000marginal}, and a doubly robust formula \citep{robins1994estimation, robins1995analysis, rotnitzky1998semiparametric, kang2007demystifying}. These formulas are functionals of measures on observed data;
    \item We deepen our understanding the influences of our identification assumptions on the set of observed data measures by showing they do not pose any restrictions and hence being nonparametric;
    \item We also provide a first-order understanding of the identification functionals on the the set of observed data measures by deriving the efficient influence curve (EIC) through semiparametric theory \citep{bickel1993efficient, van2000asymptotic, tsiatis2006semiparametric, kosorok2008introduction}.
\end{enumerate}
While this paper builds a population-level framework, it does not explore estimation or associated inferential results, which are beyond the scope of this study. With that said, our framework does pave the way for future estimation framework, by suggesting feasible estimators, possibly with efficiency. It also suggest ways to numerically approximate the causal parameter of interest, a non-trivial task for causal inference for FLD, which we illustrate through Monte-Carlo simulations. As a byproduct, our framework subsumes many existing identification methods for non-functional longitudinal data as special cases, thereby providing a unified longitudinal data causal inference framework.

The paper is organized as follows.
Section \ref{sec:pre} prepares the readers with notation and examples. Section \ref{sec:define} the potential outcome framework and define the causal parameter of interest as a functional of measures on unobserved potential outcomes, for FLD under general DTRs. Section \ref{sec:functional} presents our causal identification framework for FLD under ADTR, including the key assumptions and the derivation of the g-computation and inverse probability weighting (IPW) formulas. In this Section we also discuss the recursive representations of the g-computation and IPW processes, presents the doubly robust formula. Section \ref{sec:semi} proves that our identification assumptions do not impose any restrictions on the set of observed data distributions. Therein we also develop the semiparametric theory for our proposed estimators, deriving the efficient influence function (EIF) and discussing its implications. We conduct Monte Carlo simulations to examine our framework in Section \ref{sec:simu}, which also implies that our framework suggests a numerical approach to approximate the causal quantity of interest. In Section \ref{sec:related} we give a literature review. Finally, Section \ref{sec:dis} concludes the paper with a summary of our main contributions and a discussion of future research directions. In the appendices, Our identification assumptions are based on the so-called ``coarsened at random assumption,'' which is stronger than the classical ``sequential randomization assumption (SRA).'' Section \ref{sec:weak} provides an alternative weaker but less interpretable set of identification assumptions extending SRA, under which the results in our paper continue to hold. For readers who are less familiar with longitudinal data method or who are familiar with classical method but interested in drawing connections, we provide a review of existing methods for causal inference mathematically in classical longitudinal data in the supplementary material, highlighting their limitations when applied to FLD. Therein we also review existing approaches \citep{ying2022causal} designed for FLD, which yet did not accommodate DTRs. Proofs of our results are also given in the supplementary material.

\section{Preliminaries}\label{sec:pre}
Throughout the rest of the paper, suppose that there is a longitudinal study during 0 to $\infty$. We first define stochastic processes:
\begin{enumerate}
    \item $A(t)$ are treatments received at time $t$, which can be multi-dimensional binary, categorical or functional. We write $\bar A(t) = \{A(s): 0 \leq s \leq t\}$ and abbreviate $\bar A = \bar A(\infty)$. We define $\cA$ as the set of all possible values of $\bar a$. 
    \item $L(t)$ are measured confounders at time $t$, which can be multi-dimensional binary, categorical, or functional. We write $\bar L(t) = \{L(s): 0 \leq s \leq t\}$ and abbreviate $\bar L = \bar L(\infty)$. For convenience, we write $\bar Y$, a subset of $\bar L$, as the outcomes of interest. We include the outcomes into covariates to reflect the fact that both past outcomes and covariates influence current treatment, which can affect both future outcomes and covariates.
    \item To accommodate mortality and censoring, typically encountered in longitudinal studies, we consider $T$ as a time-to-event endpoint, for instance, death, and $C$ as the right censoring time, both of which are terminating states. Define $X = \min(T, C)$ as the censored event time. To ease notation, we include $T$ into $\bar Y$ and set $C = \infty$ whenever $T \leq C$.
    \item Note that $A(t)$ and $L(t)$ are not observed for $t \leq X$ or defined for $t \leq T$. We can offset $A(t) = A(X)$ and $L(t) = L(X)$ for observed data whenever $t > X$ and $A(t) = A(T)$ and $L(t) = L(T)$ for censoring free data whenever $t > T$. Therefore, $A(t)$ and $L(t)$ is always well defined and observed for any $t \geq 0$. Also, the information of $X$ is absorbed into $A(t)$ and $L(t)$.
    \item We introduce the counterfactual covariates $L_{\bar a}(t)$ and counterfactual outcomes $Y_{\bar a}(t)$, for any $\bar a \in \cA$, as the covariates and outcomes if the treatment were to set as $\bar a$. These counterfactuals are defined for deterministic and static treatment regimes $\bar a$ (see Section \ref{sec:TR}) and considered to be defined at baseline. We write $\bar L_{\cA} = \{\bar L_{\bar a}\}_{\bar a \in \cA}$ and $\bar Y_{\cA} = \{\bar Y_{\bar a}\}_{\bar a \in \cA}$.
    \item The observed data are $(C, \bar A, \bar L) = (C, \bar A(X), \bar L(X))$ and the uncensored full data are $(C, \bar A(T_{\cA}), \bar L_{\cA}(T_{\cA}))$.
    \item Define $\sF_t = \sigma\{\bar A(t), \bar L(t)\}$ as the natural filtration and $\sF_{t-} = \sigma(\cup_{0 \leq s < t}\sF_s)$. We also define $\sG_t = \sigma\{A(t), \sF_{t-}\}$ as a one-step treatment aware filtration. Note that $X$ is a stopping time of $\sF_t$ and $\sG_t$, thus we may write $\sF_X$ and $\sG_X$. We manually write $\sF_{0-}$ and $\sG_{0-}$ as the trivial sigma algebra for convenience. 
    \item We use the upper case for random variables and the lower case for their realized values. 
\end{enumerate}
We then introduce notation for measures:
\begin{enumerate}
    \item When there is no confusion, we may write $\P(\textup{d}\bar a\textup{d}\bar l)$ to represent the distribution on the path space induced by the stochastic processes and the probability measure $\P$ on the sample space $\Omega$. Note that this is not a density. This notation is well adopted by probabilists \citep{bhattacharya2007basic, durrett2019probability} and also statisticians \citep{gill2001causal}. 
    \item We assume the event space is Polish so that conditional probability can be chosen to be regular. We understand conditional distribution as a function over a sigma algebra multiplied with a path set. For instance, $\P(\textup{d}\bar a|\bar l)$ can be seen as a function over Borel space generated by $\{\bar a\}$ and the path set $\{\bar l\}$. Importantly, conditional distribution is defined almost surely and one needs to take extra caution when replacing and intervening treatment distributions when conducting causal inference.
    \item For any two probability distribution $\P$ and $\P'$ on $\Omega, \sF$, and two temporary random variables $U, U'$, we write $\frac{\textup{d}\P'}{\textup{d}\P}(u'|u)$ as the Radon-Nikodym derivative. It is understood as a function of $(u', u)$ such that for any measurable function $f(u', u)$, $\int f(u', u)\frac{\textup{d}\P'}{\textup{d}\P}(u'|u)\frac{\textup{d}\P'}{\textup{d}\P}(u)\P(\textup{d}u'\textup{d}u) = \int f(u', u)\P'(\textup{d}u'\textup{d}u)$.
    \item Note that conditional probability is only uniquely defined almost surely. Therefore throughout this draft, otherwise stated, for the measure zero subset where the conditional probability is not uniquely defined, we set the conditional probability to be zero.
    \item We use $\|\cdot\|_\textup{TV}$ to represent the total variation norm over the space of signed measures of the path space, which is a Banach space.
    \item We define G\^ateaux derivative. For any two probability distribution $\P$ and $\P'$ on $(\Omega, \sF)$, we write $\frac{\partial \P'}{\partial \P}$ as the G\^ateaux derivative of $\P'$ to the direction by $\P$, understood in the linear space of all signed measure on $(\Omega, \sF)$.
\end{enumerate}
Finally, a partition $\Delta_K[0, \infty]$ on $[0, \infty]$ is defined as a finite sequence of $K + 1$ numbers of the form $0 = t_0 < \cdots < t_K = \infty$. The mesh $|\Delta_K[0, \infty]|$ of a partition $\Delta_K[0, \infty]$ is $\max_{i = 0, \cdots, K - 1}(t_{j + 1} - t_j)$, representing the maximum gap length of the partition.

\section{Causal Parameters of Interest}\label{sec:define}
Here we first rigorously define counterfactual variables for any DTRs. For any partition $\Delta_K[0, \infty]$, we define a hypothetical treatment $A_{\bbG_{\Delta_K[0, \infty]}}(t)$ as a random draw from
\begin{equation}
    \bar A_{\bbG_{\Delta_K[0, \infty]}}(t_0) \sim \bbG\{\textup{d}a(t_0)\},
\end{equation}
we then define $Y_{\bbG, \Delta_K[0, \infty]}(t_0)$ as plugging in $\bar A_{\bbG_{\Delta_K[0, \infty]}}(t_0)$ into $Y_{\bar a(t_0)}(t_0)$ as
\begin{equation}
    Y_{\bbG, \Delta_K[0, \infty]}(t_0) = Y_{\bar A_{\bbG_{\Delta_K[0, \infty]}}(t_0)}(t_0).
\end{equation}
Therefore we can iteratively define
\begin{equation}
    \bar A_{\bbG_{\Delta_K[0, \infty]}}(t_1) \sim \bbG\{\textup{d}\bar a(t_1)|\bar A(t_0), \bar L(t_0), Y_{\bbG, \Delta_K[0, \infty]}(t_0)\},
\end{equation}
\begin{equation}
    Y_{\bbG, \Delta_K[0, \infty]}(t_1) = Y_{\bar A_{\bbG_{\Delta_K[0, \infty]}}(t_1)}(t_1).
\end{equation}
We iterate until infinity, so we have $\bar Y_{\bbG, \Delta_K[0, \infty]}$ defined. Now, we can let $|\Delta_K[0, \infty]| \to 0$ and define
\begin{Def}[Counterfactuals under DTRs]\label{def:counterfactualDTRs}
When the distributions $Y_{\bbG, \Delta_K[0, \infty]}$ converge to the same limit in the total variation sense, we define a realization of the limit distribution as $\bar Y_{\bbG}$.
\end{Def}
The parameter of interest is defined as
\begin{equation}\label{eq:targetparameter}
    \E\{\nu(\bar Y_{\bbG})\},
\end{equation} 
where $\nu$ is a user-specified bounded continuous function. This estimand includes those considering the marginal mean of outcomes under a STR in \citet{ying2024causality}.

Now we investigate when $\bar Y_{\bbG}$ is well-defined. In fact, when 
\begin{enumerate}
    \item either there exists a bounded function $\eps(t, \eta) > 0$ with $\int_0^\tau \eps(t, \eta)dt \to 0$ as $\eta \to 0$, such that $t \in [0, \infty]$, $\eta > 0$, 
\begin{equation}
    \E\left\{\|\bbG_1(\textup{d}\bar a|\sF_{t + \eta}) - \bbG_1(\textup{d}\bar a|\sF_{t})\|_{\textup{TV}}\right\} \leq \eps(t, \eta).
\end{equation}
    \item or $\bbG_2$ is discrete.
\end{enumerate}
Like for Riemannian integral, condition 1 or 2 serves as the corresponding ``absolute contuity'' or ``discrete'' for a function to be ``Riemannian integrable''. Condition 1 typically holds for continuous-path stochastic process like Wiener measure while condition 2 holds for point measure like Poisson process.
\begin{prp}[Existence of counterfactual outcomes under DTR for FLD]\label{prp:dtrdef}
Under either conditions listed above, the distributions of $\bar Y_{\bbG, \Delta_K[0, \infty]}$ converge in the total variation sense to the same point whenever $|\Delta_K[0, \infty]| \to 0$.
\end{prp}

We list some examples of the intervention $\bbG$ borrowed from \citet{ying2024causality} below:
\begin{itemize}
    \item When the causal outcome under a specific regime $\bar a$ is of interest, for instance, all patient were under treatment, the point mass (delta) measure $\mathbbm{1}(\bar A = \bar 1)$ can be considered;
    \item Though the data are allowed to be functional and the underlying data generating mechanism can have uncountably infinite number of treatment-confounder feedbacks, a finite-dimensional distribution intervention can still be considered, for example, intervening dosage of certain drug hourly or daily;
    \item Likewise, a marked point process measure represents intervening both dosage and frequency of usage for certain drugs;
    \item If considering certain fluid intake that is continuously used, one might leverage stationary process measure that allows noise of fluid usage yet conforms to time regularity;
    \item One may also consider continuous Gaussian process (including Wiener measure, also known as Brownian motion) as a typical example considered in stochastic processes.
\end{itemize}
Note that all interventions now can depend on historical covariates and treatment, and also the observed treatment regimes.

\section{Identification}\label{sec:functional}
\subsection{Identification assumptions}
In this section, we provide identification formulas for \eqref{eq:targetparameter}, that is, functional of the observed data distribution $\P(\textup{d}a\textup{d}l)$ that equals \eqref{eq:targetparameter}.
For any sequences of partitions $\{\Delta_K[0, \infty]\}_{K = 1}^\infty$ with $|\Delta_K[0, \infty]| \to 0$ as $K \to \infty$, we have the following decomposition
\begin{align}
    \P(\textup{d}c\textup{d}\bar a \textup{d}\bar l) &= \prod_{j = 0}^{K - 1}
    \P\{C \leq t_{j + 1}|\bar a(t_{j + 1}), \sF_{t_j}\}^{\mathbbm{1}(C \leq t_{j + 1})}\\
    &\P\{C > t_{j + 1}|\bar a(t_{j + 1}), \sF_{t_j}\}^{\mathbbm{1}(C > t_{j + 1})}\\
    &\P\{\textup{d}\bar l(t_{j + 1})|\bar a(t_{j + 1}), \sF_{t_{j}}\}\P\{\textup{d}\bar a(t_{j + 1})|\sF_{t_j}\}.
\end{align}
We define
\begin{align}
    \P_{\bbG, \Delta_K[0, \infty]}(\textup{d}c\textup{d}\bar a \textup{d}\bar l) 
    &= \mathbbm{1}(C = \infty)\prod_{j = 0}^{K - 1}\P\{\textup{d}\bar l(t_{j + 1})|\bar a(t_{j + 1}), \sF_{t_{j}}\}\bbG\{\textup{d}\bar a(t_{j + 1})|\sF_{t_{j}}\}.
\end{align}

Our first assumption is the generalized ``full conditional exchangeability'' assumption in causal inference \citep{robins2008estimation}, which generalized ``coarsening at random'' assumption in missing data literature \citep{heitjan1991ignorability}. 
\begin{assump}[Full sequential randomization]\label{assump:fullSRA}
The treatment assignment is independent of the all potential outcomes and covariates given history, in the sense that there exists a bounded function $\eps(t, \eta) > 0$ with $\int_0^\infty \eps(t, \eta)dt \to 0$ as $\eta \to 0$, such that for any $t \in [0, \infty]$, $\eta > 0$, 
\begin{equation}
    \E(\|\P\{\textup{d}\bar l_{\cA}|\bar A(t + \eta), \sF_t\} - \P\{\textup{d}\bar l_{\cA}|\sF_t\}\|_{\textup{TV}}) < \eps(t, \eta).
\end{equation}
\end{assump}

This assumption requires that the treatment assignment mechanism at any given time point is independent of all future potential outcomes and covariates, given the observed history up to that time point. To assist Assumption \ref{assump:fullSRA}, we assume:
\begin{assump}[Full consistency]\label{assump:fullcons}
\begin{equation}
    \bar L = \bar L_{\bar A}.
\end{equation}
\end{assump}
It states that, the observed covariate value under any TR is identically distributed with the counterfactual covariate value under the same TR.

We also have the positivity Assumption: 
\begin{assump}[Full positivity]\label{assump:fullpos}
\begin{equation}
\P(\textup{d}x_{\bar a}\textup{d}l_{\bar a})\bbG(\textup{d}\bar a|\bar l_{\cA})\delta_{\bar a} \ll \P(\textup{d}x\textup{d}\delta\textup{d}\bar a\textup{d}\bar l),
\end{equation}
almost surely.
\end{assump}
Furthermore, we need an assumption over the censoring mechanism to eliminate the censoring bias. We consider the conditionally independent censoring assumption in \citet{van2003unified, rotnitzky2005inverse, tsiatis2006semiparametric, ying2022causal}. Define the full data censoring time hazard as
\begin{equation}
    \lambda_C\{t|T, \overline A(T), \overline L(T)\} = \frac{\partial}{\partial \textup{d}t}\P\{C \leq t + \textup{d}t|C > t, T, \overline A(T), \overline L(T)\}/\textup{d}t.
\end{equation}
The following conditional independent censoring assumption requires that the full data censoring time hazard at time $t$ only depends on the observed data up to time $t$.
\begin{assump}[Conditional independent censoring]\label{assump:CIR}
The censoring mechanism is said to be conditionally independent if
\begin{equation}\label{eq:car}
    \lambda_C\{t|T, \overline A(T), \overline L(T)\} = \lambda_C\{t|\sF_t\}\mathbbm{1}(T > t).
\end{equation}
\end{assump}
The conditional independent censoring assumption requires that the censoring mechanism depends only on the observed data history and not on any future potential outcomes or covariates. This assumption ensures that the censoring process does not introduce bias into the estimation of causal effects, as it is independent of the potential outcomes given the observed data.

\subsection{Identification formulas}
\begin{Def}[G-computation process]\label{def:gprocess}
Under Assumptions \ref{assump:fullSRA} and \ref{assump:fullpos}, define
\begin{equation}
    H_\bbG(t) = H_\bbG(t; \sG_t) = \E_\bbG\{\nu(\bar Y)|\sG_t\},
\end{equation}
as a projection process and apparently a $\P_\bbG$-martingale. We call $H_\bbG(t)$ the {\it g-computation process}. Note that
\begin{equation}
    H_\bbG(0-) = \E_\bbG\{\nu(\bar Y)\}.
\end{equation}
\end{Def}

\begin{thm}[G-computation formula]\label{thm:gformulapoint}
Under Assumptions \ref{assump:causallyintervenable} - \ref{assump:CIR}, the marginal mean of transformed potential outcomes under a user-specified regime $\bbG$ \eqref{eq:targetparameter} is identified via a g-computation formula as
\begin{equation}
    \E\{\nu(\bar Y_{\bbG})\} = H_\bbG(0-).
\end{equation}
\end{thm}
This theorem provides a formal identification result for the marginal mean of transformed potential outcomes under a user-specified dynamic treatment regime. The g-computation formula expresses this causal quantity in terms of the observed data distribution and the target treatment regime, enabling its estimation from the available data. The proof of this formula follows a similar structure to that presented in \citet{ying2022causal}, with additional considerations for the dynamic treatment regime setting.

\begin{Def}[Inverse probability weighting process]\label{def:IPWprocess}
Under Assumptions \ref{assump:fullSRA} and \ref{assump:fullpos}, define
\begin{equation}
    Q_\bbG(t) = Q_\bbG(t; \sG_t) = \left.\E\left(\frac{\textup{d}\P_\bbG}{\textup{d}\P}\right|\sG_t\right),
\end{equation}
as a Radon-Nikodym derivative at any time $t$ and apparently a $\P$-martingale. We call $Q_\bbG(t)$ the {\it inverse probability weighting process}. Note that
\begin{equation}
    Q_\bbG(X) = \left.\E\left(\frac{\textup{d}\P_\bbG}{\textup{d}\P}\right|\sG_{X}\right) = \frac{\textup{d}\P_\bbG}{\textup{d}\P}, ~~~Q_\bbG(0-) = 1.
\end{equation}
\end{Def}

\begin{thm}[Inverse probability weighting formula]\label{thm:ipwpoint}
Under Assumptions \ref{assump:causallyintervenable} - \ref{assump:CIR}, the marginal mean of transformed potential outcomes under a user-specified regime $\bbG$ \eqref{eq:targetparameter} is identified via an inverse probability weighting formula as
\begin{equation}
    \E\{\nu(\bar Y_{\bbG})\} = \E\left\{Q_\bbG(X)\nu(\bar Y)\right\}.
\end{equation}
\end{thm}
This theorem presents an alternative identification result for the marginal mean of transformed potential outcomes using an inverse probability weighting (IPW) approach. The IPW formula expresses the causal quantity of interest as a weighted average of the observed outcomes, with weights determined by the ratio of the target treatment regime to the observed treatment distribution. This theorem provides another avenue for estimating causal effects in functional longitudinal data.

For any two $\sG_t$-adapted processes $H(t)$ and $Q(t)$, and a partition $\Delta_K[0, \infty] = \{0 = t_0 < \cdots < t_K = \infty\}$, we define
\begin{align*}
    \Xi_{\text{DR}, \bbG, \Delta_K[0, \infty]}(H, Q)
    &= Q(t_K)\nu(\bar Y)\\
    &- \sum_{k = 0}^K\left[Q(t_k)H(t_k) - Q(t_{k - 1})\int H(t_k)\bbG\{\textup{d}\bar a(t_k)|\sF_{t_{k - 1}}\}\right].
\end{align*}

We also define $\Xi_{\text{DR}, \bbG}(H, Q)$ as the $L^2$ weak limit of $\Xi_{\text{DR}, \bbG, \Delta_K[0, \infty]}(H, Q)$ whenever it exists. 
\begin{prp}[Double robustness]\label{prp:dr}
Under Assumptions \ref{assump:causallyintervenable} - \ref{assump:CIR}, and \ref{assump:rate}, and assuming
\begin{equation}
    (H_\bbG, Q_\bbG) \in \cM_{\textup{out}} \cap \cM_{\textup{trt}},
\end{equation}
when $\bbG$ is prespecified, $\Xi(H, Q)$ is a doubly robust for $\E\{\nu(\bar Y_{\bbG})\}$, in the sense that it remains unbiased when either $H_\bbG$ or $Q_\bbG$ is correct but not necessarily both. That is, for any $\sG_t$-adapted processes $Q(t)$ and $H(t)$ with $(H_\bbG, Q) \in \cM_{\textup{out}}$, $(H, Q_\bbG) \in \cM_{\textup{trt}}$ and $\sup_{t}\|H_\bbG(t) Q(t)\|_1 < \infty$, $\sup_{t}\|H(t)Q_\bbG(t)\|_{1} < \infty$, we have
\begin{equation}
    \E\{\nu(\bar Y_{\bbG})\} = \E\left\{\Xi_{\text{DR}, \bbG}(H_\bbG, Q)\right\} = \E\left\{\Xi_{\text{DR}, \bbG}(H, Q_\bbG)\right\}.
\end{equation}
\end{prp}

This theorem presents a doubly robust identification formula for the marginal mean of transformed potential outcomes under a DTR. The doubly robust formula combines the g-computation and IPW processes in a way that provides robustness against model misspecification. Specifically, the theorem shows that the doubly robust formula remains valid if either the g-computation process or the IPW process is correctly specified, even if the other is misspecified. This property enhances the reliability of causal effect estimates in practice.

A final identification formula is through the efficient influence curve (EIC). To give this formula we need to establish semiparametric theory and hence we postpone to the next section.

\section{Understanding of Probability Set and Functionals}\label{sec:semi}

In this section, we investigate the first order behavior of our identification functional around an observed data distribution $\P$. To that end, we leverage semiparametric theory. Semiparametric theory, as discussed in key works \citep{bickel1993efficient, van2000asymptotic, tsiatis2006semiparametric, kosorok2008introduction, kennedy2017semiparametric}, primarily deals with the local approximation of both the probability measure and its associated functional (or estimand). This theory emerges from applying principles of differential geometry to the concept of a probability measure ``manifold.'' It is important to distinguish this from semiparametric models, which describe probability distributions using both finite- and infinite-dimensional parameters. In the context of a differentiable manifold with a differentiable real-valued function (estimand) defined on it, the differential is understood as a linear approximation to the function at points along the tangent vectors. These tangent vectors represent an equivalence class of differentiable curves, defined through the equivalence relation of first-order contact between the curves. To calculate the differential, one can use any direction, or equivalently, any tangent vector, to determine the directional derivative. In the case of a probability measure ``manifold,'' these concepts take on explicit representations. Semiparametric theory characterizes the local neighborhood of a probability measure $\P$ using the Hellinger distance. Consequently, for any smooth curve, also known as a parametric submodel, that passes through $\P$, one can derive a corresponding tangent vector, or score, by employing Hadamard derivatives in conjunction with the Hellinger distance. The tangent space is defined as the completion of the linear spans of all tangent vectors, which forms a subspace of $L_0^2(\P)$. The directional derivative, acting as a continuous linear functional on the tangent vectors, possesses Riesz representers known as influence curves. The projection of these influence curves onto the tangent space is both unique and termed the efficient influence curve (EIC).

Understanding the differential representation of the estimand is crucial in practice. This is because the differential form has Riesz representers, which often suggests viable estimators in numerous scenarios that can substantially mitigate first-order bias. This aspect of semiparametric theory, due to its standalone significance and historical establishment in the literature for RLD and ILD \citep{rytgaard2022continuous}, underscores the need for a comprehensive semiparametric framework tailored for FLD. Historically, the development of semiparametric theory in the realm of FLD has been notably sparse, possibly attributed to the inherent complexities and challenges in managing FLD without a defined density function. Recognizing this gap, our work endeavors to pioneer in this field by formulating a generalized semiparametric theory for FLD. This advancement not only bridges a significant gap in the existing literature but also opens new avenues for enhanced data analysis and interpretation in studies involving FLD. The novelty of our approach lies in its ability to adapt and extend established principles of semiparametric theory to the unique characteristics and requirements of FLD, thereby offering a more robust and versatile analytical framework that is poised to transform the handling and understanding of such data.

\begin{thm}\label{thm:norestriction}
For any measure $\P$ over observed data $(C,  \bar A, \bar L)$ satisfying Assumption \ref{assump:continuity}, there exists a measure $\P^F$ over $(C, \bar A, \bar L, \bar L_{\cA})$ satisfying Assumptions \ref{assump:causallyintervenable} - \ref{assump:CIR} and inducing $\P$ over $(C, \bar A, \bar L)$.
\end{thm}

We show that the tangent space is full, as the space of all $L_0^2(\P)$ functions. This the ``next best'' result to full nonparametric property, that is, within the neighborhood of $\P$ satisfying Assumptions \ref{assump:causallyintervenable} - \ref{assump:CIR}, the other distributions satisfying Assumptions \ref{assump:causallyintervenable} - \ref{assump:CIR} are so rich that one cannot distinguish them with those who do no satisfy in an $L_0^2$ manner. This helps to establish the limiting IC of any ADTRs $\bbG$ as $|\Delta_K[0, \infty]| \to 0$, is the EIC, that is, the unique IC with shortest $L^2$ length and lies on the tangent space. A full tangent space trivially implies that any IC is EIC.

\begin{thm}\label{thm:tangentspace}
Under Assumptions \ref{assump:causallyintervenable} - \ref{assump:CIR}, the tangent space at $\P$ is $L_0^2(\P)$.
\end{thm}
This theorem characterizes the tangent space of the statistical model considered in the paper, which is a key concept in semiparametric theory. The theorem shows that, under the assumptions stated in the paper, the tangent space is fully nonparametric, meaning that it imposes no restrictions on the observed data distribution. This result highlights the generality and flexibility of the proposed causal inference framework for functional longitudinal data.

Finally we arrive at the EIC-based identification formula. Define
\begin{align*}
    &\Xi_{\bbG_{\Delta_K[0, \infty]}}\\
    &= \prod_{j = 0}^{K - 1}\frac{\textup{d}\bbG}{\textup{d}\P}\{\bar A(t_{j + 1})|\sF_{t_j}\}\nu(\bar Y)\\
    &+ \sum_{k = 0}^{K - 1}\prod_{j = 0}^{k - 1}\frac{\textup{d}\bbG}{\textup{d}\P}\{\bar A(t_{j + 1})|\sF_{t_j}\}\Bigg(\frac{\textup{d}[\frac{\partial}{\partial \P}\bbG - \bbG]}{\textup{d}\P}\{\bar A(t_{k + 1})|\sF_{t_k}\}\\
    &~~~\cdot\int \nu(\bar y)\prod_{j = k + 1}^{K - 1}\P\{\textup{d}\bar l(t_{j + 1})|\bar a(t_{j + 1}), \sF_{t_j}\}\bbG\{\textup{d}\bar a(t_{j + 1})|\sF_{t_j}\}\P\{\textup{d}\bar l(t_{k + 1})|\bar a(t_{k + 1}), \sF_{t_k}\}\\
    &~~~-\int \nu(\bar y)\prod_{j = k + 1}^{K - 1}\P\{\textup{d}\bar l(t_{j + 1})|\bar a(t_{j + 1}), \sF_{t_j}\}\bbG\{\textup{d}\bar a(t_{j + 1})|\sF_{t_j}\}\\
    &~~~~\P\{\textup{d}\bar l(t_{k + 1})|\bar a(t_{k + 1}), \sF_{t_k}\} 
    \frac{\partial}{\partial \P}\bbG\{\textup{d}\bar a(t_{k + 1})|\sF_{t_k}\}\\
    &~~~ + \int \nu(\bar y)\prod_{j = k}^{K - 1}\P\{\textup{d}\bar l(t_{j + 1})|\bar a(t_{j + 1}), \sF_{t_j}\}\bbG\{\textup{d}\bar a(t_{j + 1})|\sF_{t_j}\}\Bigg).
\end{align*}

\begin{thm}[Efficient Influence Curve]\label{thm:eic}
Under Assumptions \ref{assump:intervenable} and \ref{assump:contpospoint}, at the law when the derivative and the limit below can interchange,
\begin{equation}
    \frac{\partial}{\partial \theta}\lim_{|\Delta_K[0, \infty]| \to 0}\E_{\theta, \bbG_{\Delta_K[0, \infty]}}\{\nu(\bar Y)\} = \lim_{|\Delta_K[0, \infty]| \to 0}\frac{\partial}{\partial \theta}\E_{\theta, \bbG_{\Delta_K[0, \infty]}}\{\nu(\bar Y)\},
\end{equation}
and $\Xi_{\bbG_{\Delta_K[0, \infty]}}$ converge weakly in $L^2$, the efficient influence curve of $\E_{\bbG}\{\nu(\bar Y)\}$ is
\begin{equation}
    \Xi_{\bbG} - \E_{\bbG}\{\nu(\bar Y)\},
\end{equation}
where $\Xi_{\bbG}(H, Q)$ is the $L^2$ weak limit of $\Xi_{\bbG_{\Delta_K[0, \infty]}}(H, Q)$ as $|\Delta_K[0, \infty]| \to 0$. Furthermore under Assumptions \ref{assump:causallyintervenable} - \ref{assump:CIR}, the EIC-based formula is
\begin{equation}
    \E\{\nu(\bar Y_{\bbG})\} = \E\{\Xi_{\bbG}\}.
\end{equation}
\end{thm}

\section{Related Work}\label{sec:related}

\subsection{Causal Inference for Non-Functional Longitudinal Data}
\begin{figure}[H]
    \centering
    \includegraphics[scale = 0.4]{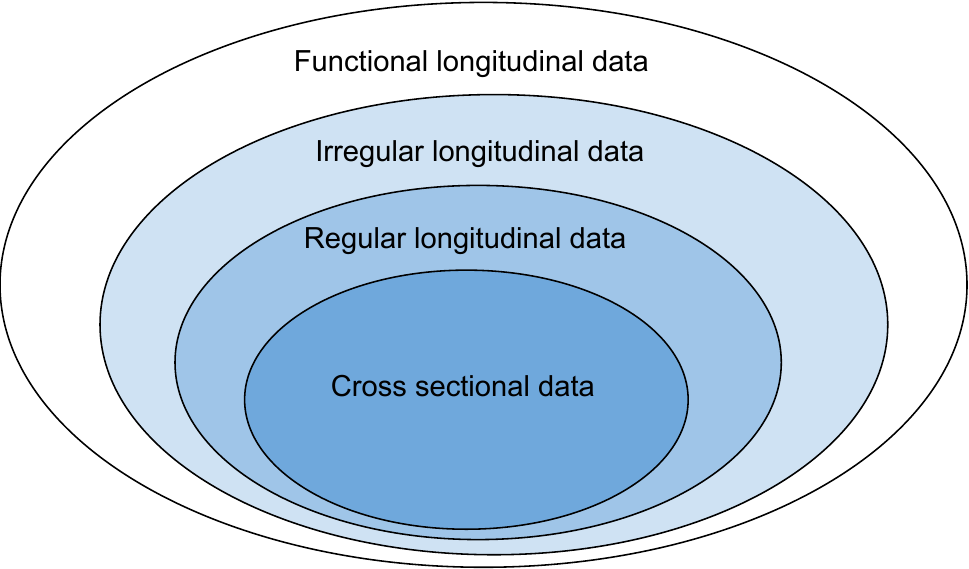}
    \includegraphics[scale = 0.55]{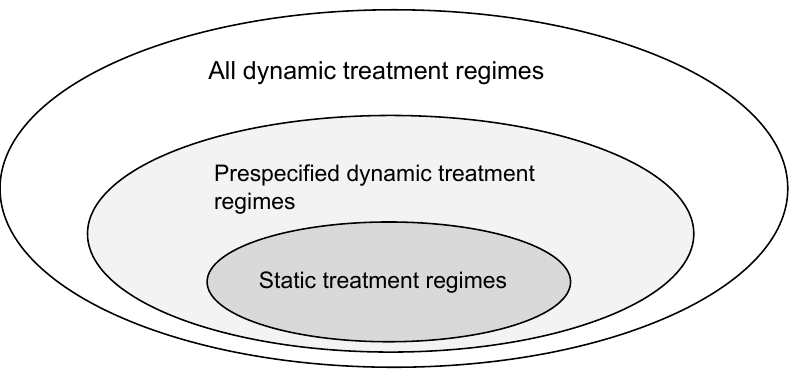}
    \caption{Left: relation between cross sectional data, regular longitudinal data, irregular longitudinal data, and functional longitudinal data. Reprinted from \citet{ying2024causality}. Right: relation between static treatment regimes, prespecified dynamic treatment regimes, and all dynamic treatment regimes. }
    \label{fig:relation}
\end{figure}

Non-functional longitudinal data can be classified into two primary types:
\begin{itemize}
    \item Regular Longitudinal Data (RLD): RLD involves data collected at set, regular intervals, though this may not accurately reflect the unpredictable nature of visit timings. 
    \item Irregular Longitudinal Data (ILD): ILD consists of data where changes in treatment and confounders occur in discrete, finite steps, resembling point processes. This type is common in pharmacoeconomic studies, where patient visits and treatment adjustments happen randomly but within finite periods. 
\end{itemize}
Figures \ref{fig:relation} provides a visual representation of relationships of different longitudinal data types.

Current causal frameworks were designed to handle RLD \citep{greenland1986identifiability, robins1986new} or ILD \citep{lok2008statistical, johnson2005semiparametric, roysland2011martingale, hu2019causal, rytgaard2022continuous}. It is evident to realize that both RLD and ILD are special cases of FLD, as shown in the left figure of Figure \ref{fig:relation}. Therefore, methods designed for FLD can be applied to RLD and ILD, but typically not vice versa. Existing frameworks tailored for RLD or ILD cannot be straightforwardly applied to FLD. Intentional or unintentional usage of existing frameworks on FLD, for instance, pruning FLD into RLD, is prone to errors and hence might lead to unwanted causal conclusions, which can have serious vital consequences and result in a waste of healthcare resources. See \citet{ying2022causal, ying2024causality} for a detailed discussion of the distinctions between the three types of longitudinal data.

\subsection{Causal Inference for Functional Data}
Existing research on causal inference has explored the realm of functional treatments and covariates within observational studies, as noted in works by \citep{miao2020average, zhang2021covariate, tan2022causal}. The data format these studies investigate is consistent with the framework of our analysis. Nonetheless, our work sets itself apart by emphasizing the temporal aspect inherent in longitudinal studies, in contrast to the primary focus on point exposure in the mentioned literature. Therefore, theoretically, by focusing on the same outcomes and setting aside the treatment-confounder feedback mechanism, our study could expand upon their findings regarding identification.

The only exceptions that investigated causal inference for functional longitudinal data are \citet{ying2022causal, ying2024causality} and \citet{sun2022causal}. However, \citet{ying2022causal, ying2024causality} only investigated static treatment regimes, without any semiparametric theory. On the other hand, \citet{sun2022causal} imposed stochastic differential equations with stringent parametric assumptions, also under static treatment regimes. This situation highlights a significant gap in methodological advancements within the field. 

\subsection{Treatment Regimes}\label{sec:TR}
A Treatment Regime (TR), also referred to as a strategy, plan, policy, or protocol, is a rule to assign treatment at each time of follow-up \citep{hernan2020causal}. 
\begin{itemize}
    \item Static vs Dynamic: A TR is called Static (STR) if it only depends on past treatment, while it is called Dynamic (DTR) if it further depends on past covariates and outcomes. 
    \item Deterministic vs Stochastic: ATR is called deterministic if it is a determined function of past information, while it is called stochastic (or random) if it is random. 
    \item Prespecified vs Actual TR dependent: 
\end{itemize}
For example, in a cross-sectional study, the average treatment effect focuses on a contrast between the averages of counterfactual outcomes if treatment were to be set to everyone or no one. This is static because it does not depend on covariates and deterministic because it is fixed. Suppose we define a rule that any patient were to receive treatment with 50\% chance, then this is a stochastic STR. If, over age 65 were not to be treated while under age 65 were to be treated, then this is a deterministic DTR. If the patient under age 65 has a 50\% chance of being treated, then this is a stochastic DTR. Apparently, deterministic TRs and STRs are special cases of stochastic TRs and DTRs, respectively.

It is important to note that previous literature almost exclusively investigates prespecified DTRs, that is, DTRs that are known and independent of the actual TR, though in most cases they have abused the name DTR. Recent works have started to consider DTRs that depend on the actual TRs. For example, \citet{haneuse2013estimation} considered ``modified treatment policy'' to answer what happens if, say, some continuous treatment dosage were to be added by some values. \citet{robins2004effects, taubman2009intervening} considered the ``threshold intervention'' to answer what happens if some continuous treatment like physical activity were to be set to a threshold if below but otherwise kept unchanged. \citet{kennedy2019nonparametric} considered the so-called ``incremental interventions'' to shift the odds of receiving treatment and answered what happens if everyone's odds of receiving treatment were, for instance, doubled compared to the actual odds. \citet{young2014identification, diaz2012population} both considered DTRs depending on actual TRs, where \citet{diaz2012population} focused on ``shift interventions'', for example, the effect of a policy that encourages people to exercise more, leading to a population where the distribution of physical activity is shifted according to certain health and socioeconomic variables. Therefore, distinguishing treatment regimes is not only a mathematical consideration but, more importantly, about what causal questions one is raising. In this paper, we investigate all DTRs (ADTRs), that is, we allow DTRs to depend on either the actual TR or not.

In Table \ref{tab:literature}, we position our paper among the literature.

\begin{table}[H]
\caption{\label{tab:literature}
A summary of literature, research gap, and our contributions.
}
\centering
\begin{tabular}{|l|l|}
\hline
\multicolumn{2}{|c|}{\multirow{2}{*}{\large Model-based methods}}  \\ 
\multicolumn{2}{|c|}{}\\\hline 
RLD & \begin{tabular}[c]{@{}l@{}}\citet{robins1997marginal, ying2023proximal};\\ \citet{vansteelandt2014structural} \end{tabular} \\ \hline
ILD &  \begin{tabular}[c]{@{}l@{}}\citet{lok2001statistical, johnson2005semiparametric, roysland2011martingale, roysland2012counterfactual};\\ \citet{hu2019causal, yang2022semiparametric, roysland2022graphical} \end{tabular}\\ \hline
FLD &  \citet{singer2008nonlinear, commenges2009general, sun2022causal}\\ \hline
\multicolumn{2}{|c|}{\multirow{2}{*}{\large Estimand-based methods for STRs}}  \\ 
\multicolumn{2}{|c|}{}\\\hline 
RLD &  \begin{tabular}[c]{@{}l@{}}\citet{rosenbaum1983central};\\ \citet{hernan2000marginal, hernan2001marginal, hernan2002estimating} \end{tabular}\\ \hline
ILD & \citet{rytgaard2022continuous}  \\ \hline
FLD & \citet{ying2022causal, ying2024causality}  \\ \hline
\multicolumn{2}{|c|}{\multirow{2}{*}{\large Estimand-based methods for PDTRs}}  \\ 
\multicolumn{2}{|c|}{}\\\hline 
RLD &  \begin{tabular}[c]{@{}l@{}}\citet{murphy2001marginal, robins2008estimation};\\ \citet{young2011comparative, chakraborty2014dynamic} \end{tabular}\\ \hline
ILD & \citet{rytgaard2022continuous}  \\ \hline
FLD & {\bf This paper} \\ \hline
\multicolumn{2}{|c|}{\multirow{2}{*}{\large Estimand-based methods for ADTRs}}  \\ 
\multicolumn{2}{|c|}{}\\\hline 
RLD &  \begin{tabular}[c]{@{}l@{}}\citet{kennedy2019nonparametric, diaz2023causal};\\ \citet{diaz2023nonparametric, wen2023intervention} \end{tabular}\\ \hline
ILD & {\bf This paper}   \\ \hline
FLD & {\bf This paper} \\ \hline
\end{tabular}
\end{table}


\section{Discussion}\label{sec:dis}
\subsection{Summary of the Paper}

By addressing these objectives and making these contributions, we aim to fill an important gap in the current literature on causal inference for complex longitudinal data. Our theoretical framework provides a foundation for the development of practical estimation methods and algorithms, paving the way for more effective treatment strategies in various domains, such as healthcare and social policy.

The identification framework and semiparametric theory can be directly applied to:
\begin{enumerate}
\item RLD with functional data observation at each regular observed time point $k$;
\item ILD with functional data observation at each jumping time $T_k$;
\item Studies on continuous time to treatment switching, initialization, or termination where timing is of interest;
\item ILD where incremental intervention, threshold intervention, or increasing intervention might be of interest;
\item Semi-Functional Longitudinal Data (SFLD) which we define as a data type standing between ILD and FLD, where the treatment process $A(t)$ is ILD and the covariate process is allowed to be FLD. 
\end{enumerate}

This paper has established a comprehensive framework for valid causal inference in continuous-time FLD with ADTRs, building upon and extending the work presented in \citet{ying2024causality, ying2022causal}. The proposed framework accommodates continuous-time progression and continuous data operation without imposing constraints on the observed data distribution. It offers sufficient assumptions for causal identification, significantly expanding upon the existing literature and encompassing scenarios with mortality and censoring. Furthermore, three distinct identification approaches have been presented: the g-computation formula, IPW formula, and an EIC formula. These contributions, along with the development of a generalized semiparametric theory for FLD, showcase the generality and nonparametric property of our framework.

\subsection{Future Directions}
Several critical areas require further investigation to advance the functional data analysis framework for longitudinal causal inference:
\begin{itemize}
    \item Assumption \ref{assump:intervenable}: Exploring alternative sufficient conditions, conducting sensitivity analyses, and examining its role in specific scenarios.
    \item Positivity Assumption: Addressing the challenges in longitudinal studies and investigating semiparametric models like marginal structural models and structural nested models to handle situations with limited data.
    \item Generalization beyond Assumption \ref{assump:contSRApoint}: Extending the framework to situations where this assumption falters, including those involving time-dependent instrumental variables and time-dependent proxies.
    \item Formal Delineation of Assumptions and Theory: Providing a more rigorous mathematical foundation for the limit in probability of pairs (H,Q) in $\cM_{\textup{out}, \bbG}$ and $\cM_{\textup{trt}, \bbG}$.
    \item Partial Identification with Discrete-Time Observations: Exploring the extent of partial identification using a discrete-time observed process for scenarios where continuous-time observations are unavailable.
\end{itemize}

By addressing these open questions and pursuing further research in these directions, we can continue to advance the field of functional data analysis for longitudinal causal inference and unlock its full potential for deriving meaningful insights from complex data scenarios. The proposed framework, along with the avenues for future research, paves the way for a more comprehensive understanding and application of causal inference techniques in the context of continuous-time longitudinal studies with ADTRs.

\begin{appendix}

\counterwithin{assump}{section}

\begin{assump}[Approximating rate]\label{assump:rate}
For any $\sG_t$-adapted process $H(t)$, and any $0 \leq s < t$,
\begin{align}
    \left|\E\left[\int H(t) \bbG\{\textup{d}\bar a(t)|\sF_s\}\P\{\textup{d}\bar l(t)|\sG_{s}\} - \E_\bbG\{H(t)|\sG_{s}\}\right]\right|
    \leq \kappa\|H(s)\|_{1}(t - s)^{\alpha},
\end{align}
for some constant $\kappa > 0$ and $\alpha > 1$.
\end{assump}

\section{Weak Assumptions}\label{sec:weak}
Throughout the paper, we consider sufficient and interpretable assumptions. In this section, we consider weaker assumptions that our framework continues to hold.

\begin{assump}[Causally Intervenable]\label{assump:causallyintervenable}
The distributions of $\bar Y_{\bbG, \Delta_K[0, \infty]}$ converge in the total variation sense to the same point whenever $|\Delta_K[0, \infty]| \to 0$, we then denote a realization of the limit as $\bar Y_{\bbG}$. 
\end{assump}
This assumption ensures that the counterfactual outcome process under a DTR $\bbG$ is well-defined, regardless of the specific partitioning of the time. It guarantees that the sequence of counterfactual outcome distributions converges to a unique limiting distribution as the mesh of the partition approaches zero. This property is crucial for making meaningful causal inferences in the FLD setting.

The following assumptions are weaker than the assumptions we gave but necessary for identifying the causal effect of interest in the context of dynamic treatment regimes:
\begin{assump}[Intervenable]\label{assump:intervenable}
The measures $\P_{\bbG, \Delta_K[0, \infty]}(\textup{d}c\textup{d}\bar a \textup{d}\bar l)$ converges to the same (signed) measure in the total variation norm on the path space, regardless of the choices of partitions, in which case we may write the limit as
\begin{equation}
    \P_{\bbG} = \P_\bbG(\textup{d}c\textup{d}\bar a \textup{d}\bar l).
\end{equation}
That is,
\begin{equation}
    \|\P_{\bbG, \Delta_K[0, \infty]} - \P_\bbG\|_{\textup{TV}} \to 0,
\end{equation}
as $|\Delta_K[0, \infty]| \to 0$. We call $\P_{\bbG}$ the target distribution.
\end{assump}
The intervenable assumption guarantees that the probability measures induced by the discretized DTR converge to a unique limiting measure, referred to as the target distribution, as the mesh of the partition approaches zero. This assumption is essential for ensuring the existence and uniqueness of the causal effect of interest in the functional longitudinal data setting.
\begin{assump}[Positivity]\label{assump:contpospoint}
The target distribution $\P_\bbG$ induced by the target regime $\bbG$ is absolutely continuous against $\P$, that is,
\begin{equation}
    \P_\bbG \ll \P,
\end{equation}
where we may write
\begin{equation}
    \frac{\textup{d}\P_\bbG}{\textup{d}\P} = \frac{\textup{d}\P_\bbG}{\textup{d}\P}(c, \bar a, \bar l)
\end{equation}
as the corresponding Radon-Nikodym derivative.
\end{assump}
The positivity assumption requires that the target distribution induced by the dynamic treatment regime is absolutely continuous with respect to the observed data distribution. In other words, for any treatment trajectory that can occur under the target regime, there must be a positive probability of observing that trajectory in the actual study. This assumption is necessary to ensure that causal effects are estimable from the observed data. The following consistency assumption is stated differently as its classic version, as the potential outcome under ADTR is defined on the distribution level.
\begin{assump}[Consistency]\label{assump:contconspoint}
\begin{equation}
    \bar Y = \bar Y_{\bar A}
\end{equation}
\end{assump}
The consistency assumption states that, for any discrete dynamic treatment regime, the counterfactual outcome process under that regime is identical in distribution to the observed outcome process when the observed TR matches the TR considered. This assumption connects the counterfactual outcomes to the observed outcomes, allowing for the estimation of causal effects from the observed data.
\begin{assump}[Sequential randomization]\label{assump:contSRApoint}
There exists a bounded function $\eps(t, \eta) > 0$ with $\int_0^\infty \eps(t, \eta)dt \to 0$ as $\eta \to 0$, such that for any DTR $\bbG$, $t \in [0, \infty]$, $\eta > 0$, 
\begin{equation}
    \left\|\E\{\nu(\bar Y_{\cA})|\sF_t\} - \E\{\nu(\bar Y_{\cA})|\bar A(t + \eta), \sF_t\}\right\|_1 < \eps(t, \eta).
\end{equation}
\end{assump}
This assumption extends the sequential randomization assumption from the discrete-time setting to the continuous-time setting. It states that, at any given time point, the treatment assignment mechanism depends only on the observed history up to that time point, and not on any future potential outcomes or covariates. This assumption is crucial for identifying causal effects in the presence of time-varying confounders in functional longitudinal data.

\subsection{Other regularity conditions}\label{sec:interpret}

\begin{prp}\label{prp:car}
Under Assumptions \ref{assump:intercontinuity} - \ref{assump:fullcons}, 
\begin{equation}
    \|\P_{\bbG, \Delta_K[0, \infty]}(\textup{d}c\textup{d}\bar a \textup{d}\bar l) - \P(\textup{d}l_{\bar a})\bbG(\textup{d}\bar a|\bar l_{\cA})\mathbbm{1}(c = \infty)\|_{\textup{TV}} \to 0,
\end{equation}
whenever $|\Delta_K[0, \infty]| \to 0$. It is immediate that Assumption \ref{assump:contpospoint} is equivalent to Assumption \ref{assump:fullpos}. Assumption \ref{assump:intervenable} holds and Assumption \ref{assump:contSRApoint} holds for any bounded function $\nu(\cdot)$.
\end{prp}
This proposition demonstrates that, under the full consistency and full conditional exchangeability assumptions (Assumptions 9 and 8), the g-computation and IPW identification results (Theorems 1 and 2) hold for the full set of counterfactual covariates and treatments. This proposition extends the identification results to a more general setting, allowing for the consideration of all potential outcomes and covariates in the causal analysis.

Consequently, one might choose Assumptions \ref{assump:CIR} - \ref{assump:fullpos} instead of Assumptions \ref{assump:causallyintervenable} - \ref{assump:CIR} as starting point if ``full conditional exchangeability'' is more intuitive.

We now proceed into investigating when the limit and derivative can interchange and when the $L^2$ weak limit $\Xi_{\bbG}$ exists in Section \ref{sec:semi}.

\begin{assump}[Uniform $L^2$ convergence]\label{assump:continuity2}
For any parametric submodel $\P_\theta$, $\Xi_{\theta, \bbG_{\Delta_K[0, \infty]}}$, defined as $\Xi_{\bbG_{\Delta_K[0, \infty]}}$ with $\P_\theta$ plugged in, converges in $L^2$ uniformly in $\theta$ as $|\Delta_K[0, \infty]| \to 0$.
\end{assump}

\begin{prp}\label{prp:weaklimit}
Under Assumptions \ref{assump:continuity2}, the regularity assumptions considered in Theorem \ref{thm:eic} holds, that is, the derivative and limit below can interchange,
\begin{equation}
    \frac{\partial}{\partial \theta}\lim_{|\Delta_K[0, \infty]| \to 0}\E_{\theta, \bbG_{\Delta_K[0, \infty]}}\{\nu(\bar Y)\} = \lim_{|\Delta_K[0, \infty]| \to 0}\frac{\partial}{\partial \theta}\E_{\theta, \bbG_{\Delta_K[0, \infty]}}\{\nu(\bar Y)\},
\end{equation}
and the $L^2$ weak limit $\Xi_{\bbG}$ of $\Xi_{\bbG, \Delta_K[0, \infty]}$ exists, as $|\Delta_K[0, \infty]| \to 0$.
\end{prp}
We also enhance our discussion on nonparametric properties. We consider the following assumption from \citet{ying2022causal}.
\begin{assump}[Observed TR Continuity]\label{assump:continuity}
There exists a bounded function $\eps(t, \eta) > 0$ with $\int_0^\tau \eps(t, \eta)dt \to 0$ as $\eta \to 0$, such that $t \in [0, \infty]$, $\eta > 0$, 
\begin{equation}
    \E\left\{\|\P(\textup{d}\bar a|\sF_{t + \eta}) - \P(\textup{d}\bar a|\sF_{t})\|_{\textup{TV}}\right\} \leq \eps(t, \eta).
\end{equation}
\end{assump}
See \citet{ying2022causal} for a complete discussion of this assumption.

\section{Recursive representation of two identification processes}

\subsection{Recursive representation of g-computation process}
Inspired by the EIC and setting $\bbG$ as if it were known, we get the following pragmatic identification formulas for $H_\bbG$ and $Q_\bbG$ and futher a doubly robust formula. For any two $\sG_t$-adapted processes $H(t)$ and $Q(t)$, and any partitions $\Delta_K[0, \infty] = \{0 = t_0 < \cdots < t_K = \infty\}$, we define
\begin{align*}
    &\Xi_{\textup{g-comp}, \bbG, \Delta_K[0, \infty]}(H, Q) \\
    &= Q(X)\{\nu(\bar Y) - H(X)\} + \sum_{j = 0}^{K - 1} Q(t_j)\left[\int H(t_{j + 1})\bbG\{\textup{d}\bar a(t_{j + 1})|\sF_{t_j}\} - H(t_j)\right].
\end{align*}
We also define $\Xi_{\textup{g-comp}, \bbG}(H, Q)$ as the $L^2$ weak limit of $\Xi_{\textup{g-comp}, \bbG, \Delta_K[0, \infty]}(H, Q)$ whenever it exists. We define
\begin{align*}
    &\cM_{\textup{out}} := \left\{(H, Q): \Xi_{\textup{g-comp}}(H, Q) \textup{ exists}\right\}.
\end{align*}
With some regularity conditions, 
\begin{equation}
    (H_\bbG, Q) \in \cM_{\textup{out}},
\end{equation}
for any locally bounded process $Q$.

\begin{thm}[Identification of the g-computation process]\label{thm:gformulaee}
Under Assumptions \ref{assump:intervenable}, \ref{assump:contpospoint}, and \ref{assump:rate}, for any $\sG_t$-adapted process $Q(t)$ with $(H_\bbG, Q) \in \cM_{\textup{out}}$ and $\sup_{t}\|H_\bbG(t)Q(t)\|_1 < \infty$, $\Xi_{\textup{g-comp}, \bbG}(H_\bbG, Q)$ is unbiased for zero, that is,
\begin{equation}\label{eq:trueoutee}
    \E\left\{\Xi_{\textup{g-comp}, \bbG}(H_\bbG, Q)\right\} = 0.
\end{equation}
Moreover, suppose there exists an $\sG_t$-adapted process $H(t)$ satisfying $\sup_t\|H(t)\|_{1, \bbG} < \infty$, so that for any $\sG_t$-adapted process $Q(t)$ with $\sup_{t}\|H(t)Q(t)\|_1 < \infty$, we have $(H, Q) \in \cM_{\textup{out}}$ and 
\begin{equation}\label{eq:outee}
    \E\left\{\Xi_{\textup{g-comp}, \bbG}(H, Q)\right\} = 0.
\end{equation}
Then $H(t)$ equals the g-computation process $H_\bbG(t)$ in Definition \ref{def:gprocess} for any $t \in [0, \infty]$ almost surely.
\end{thm}

This theorem establishes a set of population estimating equations that uniquely identify the g-computation process, which is a key component in the g-computation formula for estimating causal effects. The theorem shows that, under certain regularity conditions, the g-computation process is the only solution to these estimating equations. This result provides a foundation for developing estimators of the g-computation process and, consequently, the causal effects of interest.

\subsection{Recursive representation of IPW process}
For any two $\sG_t$-adapted processes $H(t)$ and $Q(t)$, and a partition $\Delta_K[0, \infty] = \{0 = t_0 < \cdots < t_K = \infty\}$, we define
\begin{align}
    &\Xi_{\textup{IPW}, \bbG, \Delta_K[0, \infty]}(H, Q) \\
    &= \sum_{j = 1}^K\left[Q(t_j)H(t_j) - Q(t_{j - 1})\int H(t_j)\bbG\{\textup{d}\bar a(t_j)|\sF_{t_{j - 1}}\}\right] \\
    &~~~+\left[Q(0)H(0) - \int H(0)\bbG\{\textup{d}\bar a(0)\}\right].
\end{align}
We also define $\Xi_{\textup{IPW}}(H, Q)$ as the $L^2$ weak limit of $\Xi_{\textup{IPW}, \Delta_K[0, \infty]}(H, Q)$ whenever it exists. We define
\begin{align*}
    &\cM_{\textup{trt}} := \left\{(H, Q): \Xi_{\textup{IPW}}(H, Q) \textup{ exists}\right\}.
\end{align*}
Likewise with some regularity conditions, we have
\begin{equation}
    (H, Q_\bbG) \in \cM_{\textup{trt}},
\end{equation}
for any locally bounded processes $H$ and $Q$. Here we provide a theorem for identifying the IPW process like Theorem \ref{thm:gformulaee}. 

\begin{thm}[Identification of the IPW process]\label{thm:IPWprocessee}
Under Assumptions \ref{assump:intervenable}, \ref{assump:contpospoint}, and \ref{assump:rate}, for any $\sG_t$-adapted process $H(t)$ with $(H, Q_\bbG) \in \cM_{\textup{trt}}$ and $\sup_{t}\|H(t)Q_\bbG(t)\|_{1} < \infty$, $\Xi_{\textup{IPW}}(H, Q_\bbG)$ is unbiased for zero, that is,
\begin{equation}\label{eq:truetrtee}
    \E\left\{\Xi_{\text{IPW}, \bbG}(H, Q_\bbG)\right\} = 0.
\end{equation}
Moreover, suppose there exists an $\sG_t$-adapted process $Q(t)$ satisfying $\sup_t\|Q(t)\|_1 < \infty$, so that for any $\sG_t$-adapted process $H(t)$ with $\sup_{t}\|H(t)Q(t)\|_{1} < \infty$, we have $(H, Q) \in \cM_{\textup{trt}}$ and 
\begin{equation}\label{eq:trtee}
    \E\left\{\Xi_{\text{IPW}, \bbG}(H, Q)\right\} = 0.
\end{equation}
Then $Q(t)$ equals the IPW process $Q_\bbG(t)$ in Definition \ref{def:IPWprocess} for any $t \in [0, \infty]$ almost surely.
\end{thm}
Similar to Theorem \ref{thm:gformulaee}, this theorem establishes a set of population estimating equations that uniquely identify the inverse probability weighting (IPW) process, which is central to the IPW formula for estimating causal effects. The theorem demonstrates that, under regularity conditions, the IPW process is the only solution to these estimating equations. This result lays the groundwork for developing estimators of the IPW process and the corresponding causal effects.

\end{appendix}

%


\bibliographystyle{agsm} 
\bibliography{ref}       



\end{document}